\documentclass[11pt]{article}
 \usepackage{epsfig}
  \usepackage{t1enc}
      \usepackage[latin1]{inputenc}
      \usepackage[english]{babel}
          \usepackage{latexsym}
  \usepackage{amssymb}
  \usepackage{euscript}
\usepackage{enumitem}

\newlist{noitemize}{itemize}{1}
\setlist[noitemize]{label={}, labelsep=0pt, leftmargin=0pt}

\begin{document}
 \setlength{\arraycolsep}{.136889em}
 \renewcommand{\theequation}{\thesection.\arabic{equation}}
\newcommand{\dto}{\stackrel{{d}}{\rightarrow}}
 \newtheorem{thm}{Theorem}[section]
 \newtheorem{propo}{Proposition}[section]
 \newtheorem{lemma}{Lemma}[section]
 \newtheorem{corollary}{Corollary}[section]
 \newtheorem{remark}{Remark}[section]
 \def\begg{\begin{equation}}
 \def\endd{\end{equation}}
 \def\ep{\varepsilon}
 \def\noo{n\to\infty}
 \def\al{\alpha}
 \def\be{\bf E}
 \def\bp{\bf P}
 \medskip
 \centerline{\Large\bf Strong Approximation of  the Anisotropic Random Walk  Revisited}
 \bigskip\bigskip
 \bigskip\bigskip
 \renewcommand{\thefootnote}{1}
 \noindent {\textbf{Endre Cs\'{a}ki} \newline
 Alfr\'ed R\'enyi Institute of Mathematics, Budapest, P.O.B. 127, H-1364, Hungary. E-mail address:
 csaki.endre@renyi.hu

 \ \bigskip
 \renewcommand{\thefootnote}{3}

  \noindent{\textbf{Ant\'{o}nia F\"{o}ldes}}
 \newline
 Department of Mathematics, College of Staten Island, CUNY, 2800
 Victory Blvd., Staten Island, New York 10314, U.S.A.  E-mail
 address: antonia.foldes@csi.cuny.edu

 \bigskip
 
 \medskip
 \noindent{\bf Abstract}\newline We study the path behavior of the
 anisotropic random walk on the two-dimensional lattice $\mathbb{Z}^2$.
Strong approximation of its components with independent oscillating Brownian motions
are proved. 
 \medskip


  \vspace{.1cm}

 \section{Introduction and main results}
 \renewcommand{\thesection}{\arabic{section}} \setcounter{equation}{0}
 \setcounter{thm}{0} \setcounter{lemma}{0}

We consider random walks on the square lattice $\mathbb{Z}^2$ of the plane
with possibly unequal symmetric horizontal and vertical step
probabilities, so that these probabilities can only depend on the value of
the vertical coordinate. In particular, if such a random walk is situated
at a site on the horizontal line $y=j\in \mathbb{Z}$, then at the next
step it moves with probability $p_j$ to either vertical neighbor, and with
probability $1/2-p_j$ to either horizontal neighbor. A substantial
motivation for studying such two-dimensional random walks on anisotropic
lattice has originated from transport problems of statistical physics.

More formally, consider the random walk
$\{{\bf C}(N)=\left(C_1(N),C_2(N)\right);\, N=0,1,2,\ldots\}$
on $\mathbb{Z}^2$ with the transition probabilities
$$ {\bf P}({\bf C}(N+1)=(k+1,j)|{\bf C}(N)=(k,j))={\bf P}
({\bf C}(N+1)=(k-1,j)|{\bf C}(N)=(k,j))=\frac{1}{2}-p_j,$$
$$ {\bf P} ({\bf C}(N+1)=(k,j+1)|{\bf C}(N)=(k,j))={\bf P}
({\bf C}(N+1)=(k,j-1)|{\bf C}(N)=(k,j))=p_j,$$
for $(k,j)\in{\mathbb Z^2}$, $N=0,1,2,\ldots$ We assume throughout the
paper that $0<p_j\leq 1/2$ and $\min_{j\in\mathbb Z} p_j<1/2$. Unless
otherwise stated we assume also that ${\mathbf C}(0)=(0,0)$.
In this paper we will have the following condition 

\begg n^{-1}\sum_{j=1}^n p_j^{-1}=2\gamma_1+o(n^{-\tau}), \qquad
n^{-1}\sum_{j=1}^n p_{-j}^{-1}=2\gamma_2+o(n^{-\tau}) \label{heyde0}\endd
as $n\to\infty$ for some constants $1< max(\gamma_1, \gamma_2)< \infty$,  and
$1/2<\tau \leq1 .$
We will point out how this condition is different from the previous  similar results.
The case $p_j=1/4,\, j=0,\pm 1,\pm 2,\ldots$ corresponds to simple
symmetric random walk on the plane. For this case we refer to Erd\H os and
Taylor \cite{ET}, 
Dvoretzky and Erd\H os 
\cite{DE},  and R\'ev\'esz \cite{RE}.
The case $p_j=1/2$ for some $j$ means that the horizontal line $y=j$ is
missing. If all $p_j=1/2$, then the random walk takes place on the $y$ axis,
so it is only a one-dimensional random walk, and this case is excluded from
the present investigations. The case however when $p_j=1/2,\,
j=\pm 1,\pm 2,\ldots$ but $p_0=1/4$ is an interesting one which is the
so-called random walk on the two-dimensional comb. In this case  $\gamma_1=\gamma_2=1.$ For this model we may
refer to Weiss and Havlin \cite{WH}, Bertacchi and Zucca \cite{BZ},
Bertacchi \cite{BE}, Cs\'aki {\it et al.} \cite{CCFR08}. In the comb model the scaling 
of the horizontal and vertical coordinates are different, namely for the first coordinate 
the scaling is of order $N^{1/4},$ so it is a so called sub-diffusion,  and can be
 approximated with an iterated Wiener process while the second coordinate
is of order  $N^{1/2},$ and hence it can be approximated with a Wiener process.

 In our paper  Cs\'aki {\it et al.} \cite{CCFR13} we considered the  case when
   both coordinates are of order $N^{1/2},$ hence can be approximated simultaneously with independent  Wiener processes.
    More precisely we investigated the case  when  in   (\ref{heyde0})  we have  $\lambda_1=\lambda_2>1.$ 
 In  Cs\'aki {\it et al.} \cite{CCFR13} we proved that

\medskip\noindent 
 {\bf  Theorem A} {\it  Under the condition {\rm (\ref{heyde0})} with  $\lambda=\lambda_1=\lambda_2>1,$   and $1/2<\tau\leq 1$,
on an appropriate probability space for the random walk
$$\{{\bf
C}(N)=(C_1(N),C_2(N)); N=0,1,2,\ldots\}$$
one can construct two independent
standard Wiener processes $\{W_1(t);\, t\geq 0\}$, $\{W_2(t);\,
t\geq 0\}$ so that, as $N\to\infty$, we have with any}
$\varepsilon>0$
\begg
\left|C_1(N)-W_1\left(\frac{\gamma-1}{\gamma}\,N\right)\right|+
\left|C_2(N)-W_2\left(\frac{1}{\gamma}\,N\right)\right|
=O(N^{5/8-\tau/4+\varepsilon})\quad {a.s.}
\endd

\medskip\noindent

    The case $\lambda_1=\lambda_2>1$ has a considerable 
complex history, we refer to the interested reader to \cite{CCFR13}. Here we just mention a few names;
Silver {\it et al.} \cite{SSL}, Seshadri {\it et al.} \cite{SLS},  Shuler
\cite{SH}, Westcott \cite{WE}. Some of the most important contribution to this topic is due to Heyde  \cite{H}, \cite{H93} and Heyde {\it et al.} \cite{HWW}. 
 Let $\{Y(t), t \geq 0\}$ be a
diffusion process on the same probability  space  as $\{C_2(n)\} $ whose 
distribution is defined by  
 $$Y(t)=W(A^{-1}(t)), \quad  t\geq 0,$$
where $\{W(t), t\geq 0 \}$ is a standard Brownian motion, $($or standard Wiener process$) $ and
 $$A(t)=\int_0^t\sigma^{-2}(W(s))\,ds$$
and
 \begin{displaymath}
\sigma^2(y)=\left\{ 
\begin{array}{ll}
& \frac{1}{\gamma_1}\, \quad {\rm for}\quad y\geq 0, \\
\\
& \frac{1}{\gamma_2}\, \quad {\rm for}\quad y<0.
\end{array}\right.
\end{displaymath}

Here $A^{-1}(\cdot)$ is the inverse of $A(\cdot)$. The process $Y(t)$ is 
called oscillating Brownian motion if $\gamma_1\neq \gamma_2$, that is a 
diffusion with speed measure $m(dy)= 2\sigma^{-2}(y)dy.$

\medskip\noindent
{\bf Remark 1.1} Observe that $A(t)$ above  is equal to
\begg 
A(t)=\gamma_1\int_0^t I(W(s)\geq 0)\,ds+\gamma_2\int_0^t I(W(s)< 0)\,ds.
\label{at}
\endd 

Let 
\begg
k^{-1}\sum_{j=1}^k p_{j}^{-1}=2\gamma_1+\ep_k \quad
k^{-1} \sum_{j=-k}^{-1} p_j^{-1}=2\gamma_2+\ep^*_k    
\label{H1}\endd
then the main result of  Heyde {\it et al.} \cite{HWW}  is that

\noindent
{\bf Theorem B} (\cite{HWW}) {\it Suppose that in {\rm (\ref{H1})} 
$\ep_k$ and $\ep^*_k$ are $o(1)$ as $k\to\infty.$ Then} 
$$
\sup_{0\leq t\leq N} | N^{-1/2}C_2({[Nt]})-Y(t)| \to 0  \quad a.s.
$$ 
Observe that here $\gamma_1$ and $\gamma_2$ might be different,  and the convergence rates are much less restrictive, but the approximation is only for the second component.

 Lets define  an arbitrary set $B\subset \mathbb Z $ such that for  
 
  \begg p_i=1/4 \quad {\rm  if} \quad  i \in B \quad  {\rm and } \quad   p_i=1/2  \quad   i  \in \mathbb Z  \smallsetminus B. \label{setb} \endd

Thus we remove from the two-dimensional integer lattice 
all the horizontal edges which do not belong to the $i$-levels in $B.$
In our paper \cite{CCFR12}  we investigated a simple random walk on the half-plane 
half-comb (HPHC) structure, which is another interesting special case  where we define  the   set
$B=\{ i=0,1,2,...\}$, that is to say, all horizontal lines under the $x$-axis 
are deleted. Our main result there reads as follows.

\medskip\noindent
{\bf Theorem C} (\cite{CCFR12}) {\it On an appropriate probability space for 
the HPHC random walk \newline $\{{\bf C}(N)=(C_1(N),C_2(N)); N=0,1,2,\ldots\}$ 
with $p_j=1/4,\, j=0,1,2,\ldots$, $p_j=1/2,\, j=-1,-2,\ldots$ one can 
construct two independent standard Wiener processes 
$\{W_1(t);\, t\geq 0\}$, $\{W_2(t);\, t\geq 0\}$ such that, as $N\to\infty$, 
we have with any} $\varepsilon>0$
$$
|C_1(N)-W_1(N-A_2^{-1}(N))|+|C_2(N)-W_2((A_2^{-1}(N))|
=O(N^{3/8+\varepsilon})\quad a.s.,
$$
{\it where} $A_2(t)=2\int_0^t I(W_2(s)\geq 0)\,ds+\int_0^t I(W_2(s)< 0)\,ds.$

\medskip\noindent

Clearly in this case $\lambda_1=2$ and $\lambda_2=1,$ $\tau$ can be selected to be 1.
In our paper Cs\'aki and F\"oldes \cite{CF21} we considered the case when the set $B$ is much more general than in Theorem C.  
Our main result in that paper can  be formulated as follows:

\medskip\noindent
{\bf Theorem D} \cite{CF21} {\it Let   \begg |B_n|:=|B\cap \{ -n, n\}|\sim c n \label{beta} \endd  with some constant $c>0,$ where  $B$ is defined in} (\ref{setb}) { \it and $|B_n|$ stands for the 
(finite) number of elements in the set} $B_n$. {\it Under the conditions {\rm (\ref{heyde0})}  with 
$\max(\gamma_1,\gamma_2)>1$
on an appropriate probability space for the random walk $\{{\bf
C}(N)=(C_1(N),C_2(N));\, \, N=0,1,2,\ldots\}$ one can construct two independent
standard Wiener processes $\{W_1(t);\, t\geq 0\}$, $\{W_2(t);\,
t\geq 0\}$ so that, as $N\to\infty$, we have with any
$\varepsilon>0$
\begg
\left|C_1(N)-W_1\left(N-A_2^{-1}(N)\right)\right|+
\left|C_2(N)-W_2\left(A_2^{-1}(N)\right)\right|
=O(N^{5/8-\tau/4+\varepsilon})\quad a.s. \label{hat}
\endd
where $A_2^{-1}(\cdot)$ is the inverse of }$A_2(\cdot).$

So in Theorem D we have exactly our condition  (\ref {heyde0}), but all  the $p_i$-s has to be 1/2 or 1/4.
 Our goal in this paper is to get a common generalization of  Theorems A,B,C and D, namely 
we only need condition  (\ref {heyde0}) and no other restrictions for the $p_i$-s.

\begin{thm} Under the conditions {\rm (\ref{heyde0})}, and   
$1\leq \gamma_2 < \gamma_1$
on an appropriate probability space for the random walk $\{{\bf
C}(N)=(C_1(N),C_2(N));\, \, N=0,1,2,\ldots\}$ one can construct two independent
standard Wiener processes $\{W_1(t);\, t\geq 0\}$, $\{W_2(t);\,
t\geq 0\}$ so that, as $N\to\infty$, we have with any
$\varepsilon>0$
\begg
\left|C_1(N)-W_1\left(N-A_2^{-1}(N)\right)\right|+
\left|C_2(N)-W_2\left(A_2^{-1}(N)\right)\right|
=O(N^{5/8-\tau/4+\varepsilon})\quad a.s. \label{hat}
\endd
{\it where} $A_2(t)=\lambda_1\int_0^t I(W_2(s)\geq 0)\,ds+\lambda_2\int_0^t I(W_2(s)< 0)\,ds.$

 and $A_2^{-1}(\cdot)$ is the inverse of $A_2(\cdot).$
\end{thm}

\noindent
{\bf Remark 1.2} If $ \gamma_1=\gamma_2>1$ then  $A_2(t)=\gamma_1t$ and our theorem  coincides with Theorem A. So we made the supposition that $1\leq \gamma_2 < \gamma_1,$  instead of  $1< max(\gamma_1, \gamma_2)< \infty,$ even though it is not necessary but makes the flow of argument easier.

\section{Preliminaries}
\renewcommand{\thesection}{\arabic{section}} \setcounter{equation}{0}
\setcounter{thm}{0} \setcounter{lemma}{0}

First we are to redefine our random walk $\{{\mathbf C}(N);\,
N=0,1,2,\ldots\}$. It will be seen that the process described right below
is equivalent to that given in the Introduction (cf. (\ref{equ}) below).

To begin with, on a suitable probability space  consider two independent
simple symmetric (one-dimensional) random walks $S_1(\cdot)$, and
$S_2(\cdot)$. We may assume that on the same probability space we have a
double array of independent geometric random variables $\{G_i^{(j)},\,
i\geq 1, j\in \mathbb{Z}\}$ which are independent from $S_1(\cdot)$, and
$S_2(\cdot),$ where $G_i^{(j)}$ has the following geometric distribution
\begg
\mathbf{P}(G_i^{(j)}=k)=2p_j(1-2p_j)^k,\,\, k=0,1,2,\ldots \label{geo}
\endd

We now construct our walk $\mathbf{C}(N)$ as follows. We will take all
the horizontal steps consecutively from $S_1(\cdot)$ and all the vertical
steps consecutively from $S_2(\cdot).$ First we will take some horizontal
steps from $S_1(\cdot)$, then exactly one vertical step from $S_2(\cdot),$
then again some horizontal steps from $S_1(\cdot)$ and exactly one vertical
step from $S_2(\cdot),$ and so on. Now we explain how to get the number of
horizontal steps on each occasion.
Consider our walk starting from the origin proceeding first horizontally
$G_1^{(0)}$ steps (note that $G_1^{(0)}=0$ is possible with probability
$2p_0$), after which it takes exactly one vertical step, arriving
either to the level $1$ or $-1$, where it takes $G_1^{(1)}$ or
$G_1^{(-1)}$ horizontal steps (which might be no steps at all) before
proceeding with another vertical step. If this step carries the walk to
the level $j$, then it will take $G_1^{(j)}$ horizontal steps, if this is
the first visit to level $j,$ it takes\ $G_2^{(j)}$ horizontal
steps, if this is its second visit at level $j$  and so on. In general, 
if we finished the $k$ -th vertical step and arrived to
the level $j$ for the $i$-th time, then it will take $G_i^{(j)}$
horizontal steps.

Let now $H_N,\, V_N$ be the number of horizontal and
vertical steps, respectively from the first $N$ steps of the
just described process. Consequently, $H_N+V_N=N$, and
$$
\left\{{\bf C}(N);\, N=0,1,2,\ldots\right\}=
\left\{(C_1(N),C_2(N));\, N=0,1,2,\ldots\right\}
$$
\begg
\stackrel{{d}}{=}\left\{(S_1(H_N),S_2(V_N));\, N=0,1,2,\ldots\right\},
\label{equ}
\endd
where $\stackrel{{d}}{=}$ stands for equality in distribution.

\bigskip
We will need the following two lemmas  from Cs\'aki at al. \cite{CCFR13} 

\noindent {\bf Lemma A}
{\it Let $\{G_i^{(j)},\, i=1,2,\ldots,n_j,\, j= 1, 2,\ldots, K\}$ be
independent geometric random variables with distribution
$${\bf P}(G_i^{(j)}=k)=\alpha_j(1-\alpha_j)^k,\quad k=0,1,2,\ldots,$$
where $0<\alpha_j\leq 1$. Put
$$
B_K=\sum_{j=1}^K\sum_{i=1}^{n_j}G_i^{(j)},\quad
\sigma^2=Var B_K=\sum_{j=1}^K\frac{n_j(1-\alpha_j)}{\alpha_j^2}.
$$
Then, for
$\lambda<-\sigma^2\log(1-\alpha_j)$ for each $j\in [1,K]$, we have
\begg
{\bf P}\left(
\left|\sum_{j=1}^K\sum_{i=1}^{n_j}
\left(G_i^{(j)}-
\frac{1-\alpha_j}{\alpha_j}\right)\right|>\lambda\right)\leq
2\exp\left(-\frac{\lambda^2}{2\sigma^2}+\sum_{\ell=3}^\infty
\frac{\lambda^\ell}{\sigma^{2\ell}}\sum_{j=1}^K\frac{n_j}{\alpha_j^\ell}
\right).
\label{exp}
\endd
}
\noindent{\bf Lemma B}
{\it Assume the conditions of {\rm Lemma A} and put
$M=\sum_{j=1}^K n_j$. For $M\to\infty$ and $K\to\infty$ assume moreover
that
\begg
K=K(M)=O(M^{1/2+\delta}),\quad \max_{1\leq j\leq K}
n_j=O(M^{1/2+\delta}),
\label{km}
\endd
for all $\delta>0$,
\begg
\frac{1}{\alpha_j}\leq c_1|j|^{1-\tau},\quad j=0, 1, 2,\ldots
\label{alpha}
\endd
for some $1/2<\tau\leq 1$ and $c_1>0$,
\begg
\sum_{j=1}^K\frac{1}{\alpha_j}=O(K),\qquad \frac{1}{\sigma}\leq
\frac{c_2}{M^{1/2}}
\label{sigma}
\endd
for some $c_2>0$. Then we have as $K,M\to\infty$,
\begg
\sum_{j=1}^K\sum_{i=1}^{n_j}G_i^{(j)}=\sum_{j=1}^K
n_j\frac{1-\alpha_j}{\alpha_j}+O(M^{3/4-\tau/4+\varepsilon})
\quad a.s.
\label{main}
\endd
for some $\varepsilon>0$.
}

\bigskip
 \noindent
{\bf Remark 2.1} Lemmas A and B in  \cite{CCFR13} are formulated  for  $-K\leq j\leq K,$ but it is obvious that the present  one can be stated and proved word  by word as it is  in  \cite{CCFR13}. 

\medskip\noindent
	Let $ \{X_i\}$ be a sequence of independent i.i.d. random variable, with $ {\bf P}(X_i=\pm 1)=1/2.$
Let $S(n)=\sum_{i=1}^n X_i.$ Then $S(n)$ is a simple symmetric walk on the line.  Its local time is defined by  $\xi(j,n)=\# \{k: 0<k<n,S(k)=j\}, \quad n=1,2,....$\,  for any integer $j.$  
Define $M_n=\max_{0\leq k \leq n}|S(k)|.$  Then we have the usual law of the
iterated logarithm (LIL) and Chung's LIL \cite{CH}.

\medskip\noindent
{\bf Lemma C}
$$\limsup_{\noo} \frac{M_n}{(2n \log\log n)^{1/2}}=1, \qquad
\liminf_{\noo}\left(\frac{\log\log
n}{n}\right)^{1/2}M_n=\frac{\pi}{\sqrt{8}}\qquad a.s. $$
\bigskip

The following result is from Heyde \cite{H}, (see also in \cite{CR83} Lemma 5)

\medskip\noindent
{\bf Lemma D}  {\it For the simple symmetric random walk for any 
$\varepsilon>0$ we have}
$$\lim_{n\to \infty}
\frac{\sup_{x\in \mathbb{Z}}|\xi(x+1,n)-\xi(x,n)|}{n^{1/4+\varepsilon}}=0\quad
a.s.
$$

For the next Lemma see Kesten \cite{K}.

\medskip\noindent
{\bf Lemma E}   {\it For the maximal local time

$$\xi(n)=\sup_{x \in \mathbb{Z}}\, \xi(x,n) $$
 we
have}
$$ \limsup_{n\to\infty} \frac{\xi(n)}{(2n\log \log n)^{1/2}}=1\quad {a.s.}$$

We need a simple lemma about the properties of $A(t)$ from Cs\'aki and   F\"oldes \cite{CF21}

\medskip\noindent
{\bf Lemma F} Consider $A(t)$ defined by (\ref{at}) and let $\alpha(t)=A(t)-t$.
 { \it Let $\gamma_1>\gamma_2\geq1.$ Then 
\begin{itemize}
\item
$A(t)  \quad {\rm and}\quad \alpha(t) $ {\rm are nondecreasing}
\item
$\gamma_2 t \leq A(t) \leq \gamma_1t   \quad {\rm and} \quad 
\frac{t}{\gamma_1}\leq A^{-1}(t) \leq \frac{t}{\gamma_2}.$ 
\end{itemize}
}

We will need the famous KMT strong
invariance principle (cf. Koml\'os {\it et al.} \cite{KMT}).

\medskip\noindent
{\bf Lemma G} 
{\it On an appropriate probability space one can construct
$\{S(n),\, n=1,2,\ldots\}$, a simple symmetric random walk on the
line and a standard Wiener process $\{W(t),\, t\geq 0\}$ such that
as $n\to\infty$,}
$$
S(n)-W(n)=O(\log n)\qquad a.s.
$$

The next lemma  is a simultaneous strong approximation result of R\'ev\'esz \cite{RE81}

\medskip\noindent
{\bf Lemma H}  {\it On an appropriate probability space for a 
simple symmetric random walk 

\noindent
$\{S(n);\, n=0,1,2,\ldots\}$ with local time
$\{\xi(x,n);\, x=0,\pm1,\pm2,\ldots;\, n=0,1,2,\ldots\}$ one can
construct a standard Wiener process $\{W(t);\, t\geq 0\}$ with local time
process $\{\eta(x,t);\, x\in\mathbb R; \, t\geq 0\}$ such that, as
$n\to\infty$, we have for any $\varepsilon>0$

\begin{itemize}
 \item $|S(n)-W(n)|=O(n^{1/4+\varepsilon})\quad {a.s.}$

and
\item 
$ \sup_{x\in\mathbb Z}|\xi(x,n)-\eta(x,n)|=O(n^{1/4+\varepsilon})
\quad {a.s.},$
\end{itemize}
simultaneously.}

Concerning the increments of the Wiener process we quote the following
result from Cs\"org\H o and R\'ev\'esz (\cite{CSR79}, page 69).

\medskip\noindent
{\bf Lemma I}  {\it Let $0<a_T\leq T$ be a non-decreasing 
function of $T$. Then, as $T\to\infty$, we have}
$$
\sup_{0\leq t\leq T-a_T}\sup_{s\leq a_T}|W(t+s)-W(t)|=
O(a_T(\log(T/a_T)+\log\log T))^{1/2} \qquad a.s.
$$
The above statement is also true if $W(\cdot)$ replaced by the simple 
symmetric random walk $S(\cdot).$
\medskip\noindent

For a simple random walk with local time $\xi(\cdot,\cdot),$ let 
\begg
\widehat{A}(n)=\gamma_1 \sum_{j=0}^{\infty} \xi(j,n)+\gamma_2 
\sum_{j=1}^{\infty} \xi(-j,n).
\label{an}
\endd
 We will need the following lemma from Cs\'aki and   F\"oldes \cite{CF21}

 \medskip\noindent
{\bf Lemma J} 
{\it On a probability space as in Lemma H}
$$|\widehat{A}(n)-A(n)|=O(n^{3/4+\ep})\quad\quad a.s.$$
 where $A(.)$ is defined in (\ref{at}).
The next lemma  is using some ideas from Heyde \cite{H}.
\begin{lemma}   Let $\{S(i)\}_{i=1}^{ \infty } $ a simple symmetric random walk with local time $\xi(i,n).$ Let $\beta_j>0 \quad j=1,2, \ldots$  Suppose that
\begg n^{-1}\sum_{j=1}^n \beta_j^{-1}=\rho+o(n^{-\tau}), \qquad
\label{h} \label{heyde2}\endd
as $n\to\infty$ for some constants $0< \rho< \infty$,  and
$1/2<\tau \leq1 .$
Then 
\begg
\sum_{j=1}^N \xi(j,N) \frac{1}{\beta_j}=\rho \sum_{j=1}^{\infty} \xi(j,N) +O(N^{5/4-\tau/2+\ep}) \quad a.s.\label{szumma}
\endd
\end{lemma}

{\bf Proof}:
Introduce the notation: 
$$\frac{1}{j}\sum_{k=1}^j \frac{1}{\beta_k}=\kappa_j   \quad j=1,2, \ldots.$$

$$\sum_j \xi(j,N)\frac{1}{\beta_j}=\sum_{j=1}^\infty
\xi(j,N)(j\kappa_j-(j-1)\kappa_{j-1})$$
$$=\sum_{j=1}^\infty j \kappa_j(\xi(j,N)- \xi(j+1,N))
$$
$$=\sum_{j=1}^\infty j(\kappa_j -\rho)(\xi(j,N)-
\xi(j+1,N))+\rho \sum_{j=1}^\infty j(\xi(j,N)- \xi(j+1,N))$$

$$=\rho \sum_{j=1}^{\infty} \xi(j,N)+
\sum_{j=1}^{\infty} j(\kappa_j-\rho)(\xi(j,N)
-\xi(j+1,N))$$
Observe that from (\ref{heyde2}) we have that
$$|j(\kappa_j- \rho)|\leq c j^{1-\tau}
$$
for some $c>0.$ Now applying Lemma C for $S (\cdot)$, and Lemma D, we get
that
$$\sum_{j=1}^{\infty} j(\kappa_j-\rho)(\xi(j,N)- \xi(j+1,N))
$$
$$=O(N^{1/4+\epsilon})\sum_{j=1}^{\max_{k\leq N}|S(k)|} j^{1-\tau}=
O(N^{1/4+\epsilon})O(N^{1-\tau/2+\ep})=O(N^{5/4-\tau/2+\ep}) \qquad  a.s.,$$
where here and throughout the paper the value of $\varepsilon$ might change
from line to line.
$\Box$

\section{Proofs}
\renewcommand{\thesection}{\arabic{section}} \setcounter{equation}{0}
\setcounter{thm}{0} \setcounter{lemma}{0}

{\bf Proof of Theorem 1.1} Recall that $H_N$ and $V_N$ are the number of horizontal and vertical steps respectively of the first $N$ steps of 
 $\{{\mathbf C}(.)\}$. First we would like to approximate  $H_N$
almost surely as  $N\to \infty.$ 

Consider the sum
$$G_1^{(j)}+G_2^{(j)}+ ...+G_{\xi_2(j,V_N)}^{(j)}$$
which is the total number of horizontal steps on the level $j$, where
$\xi_2(\cdot,\cdot)$ is the local time of the walk $S_2(\cdot)$. This
statement is slightly  incorrect if $j$ happens to be the level where the
last vertical step (up to the total of $N$ steps) takes the walk. In this
case the last geometric random variable might be truncated. However the
error which might occur from this simplification will be part of the
$O(\cdot)$ term. This can be seen as follows. Let
$$H_N^*=\sum_j\sum_{i=1}^{\xi_2(j,V_N)}G_i^{(j)},
$$
where $G_i^{(j)}$ has distribution (\ref{geo}). 
Clearly
$$
H_N^* -H_N\leq\max_jG_{\xi_2(j,V_N)}^{(j)}.
$$
Here and in the sequel
$$\sum_j=\sum_{\min_{0\leq k\leq V_N}S_2(k)\leq j\leq
\max_{0\leq k\leq V_N}S_2(k)}$$
 and
$$ \max_j=\max_{\min_{0\leq k\leq V_N}S_2(k)\leq j\leq
\max_{0\leq k\leq V_N}S_2(k)}
$$

Note that  from (\ref{heyde0}) we have
$$\frac{1}{\alpha_j}=\frac{1}{2p_j}\leq  c_1|j|^{1-\tau}\quad j=\pm1, \pm2, ...$$

\begin{eqnarray}
P(\max_j G_{\xi_2(j,V_N)}^{(j)}>N^{1/2+\delta})&\leq& 
\sum_j P(G_1^j>N^{1/2+\delta})
\leq \sum_j (1-\alpha_j)^{N^{1/2+\delta}}
\\  \nonumber
&\leq& \sum_j \exp(-\alpha_j N^{1/2+\delta})
\leq\sum_j  \exp(-cj^{\tau-1}N^{1/2+\delta})
\\  \nonumber
&\leq& N^{1/2+\delta^*} \exp(-cN^{\tau-1}N^{1/2+\delta}) \\  \nonumber
&\leq& N^{1/2+\delta^*} \exp(-cN^{\tau-1/2+\delta})\leq \exp(-cN^{\ep}) \label{bound}
\end{eqnarray}
with some small $\ep>0, \delta>0, \delta^*>0.$ In the last line  we used that $1/2<\tau\leq1$ and  Lemma C.
Here and in what follows the value of $c$ can be different from line to line.
By Borel Cantelli we have now that for $N\to  \infty$ 
\begg H_N^*-H_N\leq N^{1/2+\delta} \label{haen} \endd 
almost surely for any  $\delta>0.$
 
 \bigskip
 \noindent
 The next step is to show that
$ \sum_j \xi_2(j,V_N)\frac{1}{2p_j}$ is close to $A(V_N).$

To see this we apply  Lemma 2.1 for the vertical walk  $S_2(.) $ two times (separately for positive and negative $j$ indices) with $ \beta_j= 2p_j=\alpha_j$  for  $j=\pm1, \pm2, \dots$ and   
 $\rho=\gamma_1$ and $\gamma_2$ respectively for positive and negative indices and with 
 $V_N $ instead of $N$ to conclude, that

\begin{eqnarray}\sum_ j\frac{\xi_2(j,V_N)}{\alpha_j} &=& \gamma_1 \sum_{j=1}^\infty \xi_2(j,V_N) +
\gamma_2 \sum_{j=1}^\infty \xi_2(-j,V_N) + \xi_2(0,V_N) \frac{1}{2p_0} + O(N^{5/4-\tau/2+\ep})  \\
\nonumber
&=&
\widehat{A}_2(V_N)+ O(N^{5/4-\tau/2+\ep})=A_2(V_N)+ O(N^{5/4-\tau/2+\ep}).\label{gammak}
\end{eqnarray}
where in the last line we used Lemmas E  and J and  that $V_n\leq N.$

The rest of proof will be different for $\gamma_2 > 1$ and for $\gamma_2=1.$

Consider first the case $\gamma_2 > 1.$
In this case we will apply Lemma B  twice  for 
$$H_N^*=\sum_j \left(G_1^{(j)}+G_2^{(j)}+
...+G_{\xi_2(j,V_N)}^{(j)}\right).$$

  Introduce
 $V_N(+)=\sum_{j=0}^{\infty} \xi_2 (j,V_N)$ and
$V_N(-)=\sum_{j=1}^{+\infty} \xi_2(-j,V_N)$, and   
 and let  $K=\max_{0\leq k\leq V_N}|S_2(k)|.$  Define
$$\sigma_1^2=\sum_{j=0}^K \frac{\xi_2 (j,V_N)(1-\alpha_j)}{\alpha_j^2}, \quad {\rm and} \quad
\sigma_2^2=\sum_{j=1}^K\frac{\xi_2 (-j,V_N)(1-\alpha_{-j})}{\alpha_{-j}^2}.$$
Let  $M=V_N(+)$ and $M=V_N(-) $ for the first and second application respectively
 and
$n_j=\xi_2(j, V_N)$, $\alpha_j=2p_j$, $j=0, 1, 2,\ldots$  for the first and $n_j=\xi_2(j, V_N)$, 
$\alpha_j=2p_j$, $j=-1, -2, \dots   -K$  indices for the second application.

We need to check all the assumptions of the lemma in both cases . (\ref{km}) follows from Lemma C and Lemma E, and
 (\ref{alpha}) follows from(\ref{heyde0}). 
The first part of (\ref{sigma}) follows from (\ref{h}) with $\rho$  equal $\gamma_1$  and $\gamma_2$ respectively. It remains to
verify the second part of (\ref{sigma}). 
Using Lemma 2.1   we have almost surely as  $N\to\infty$,
\begin{eqnarray}
\frac{\sigma_1^2}{V_N(+)}&=&\frac{1}{V_N(+)}\sum_{j\geq 0}\frac{\xi_2(j,V_N)(1-2p_j)}
{(2p_j)^2}=\frac{1}{V_N)(+)}\sum_{j\geq 0}\frac{\xi_2(j,V_N(+))(1-2p_j)}
{(2p_j)^2}\\ \nonumber &\geq& \frac{1}{V_N(+)}\sum_{j\geq 0 }\frac{\xi_2(j,V_N(+))(1-2p_j)}{2p_j}
 =\frac{1}{V_N(+)}\sum_{j\geq 0}\frac{\xi_2(j,V_N(+))}{2p_j}-1 \to \gamma_1-1>0.
\end{eqnarray}

 where the last inequality follows from our supposition of $\gamma_1>1.$ The corresponding argument 
 for $\sigma_2^2$  goes the same way using now that $\gamma_2>1$ as well. So we checked all the conditions of Lemma B and we can conclude that 
we have almost surely, as $N\to\infty$,
\begin{eqnarray}
H_N^*(+):=\sum_{j\geq 0} \left(G_1^{(j)}+G_2^{(j)}+\ldots
+G_{\xi_2(j,V_N)}^{(j)}\right) \label{pos}
=\sum_{j\geq 0}\xi_2(j,V_N)\frac{1-2p_j}{2p_j}+O(N^{3/4-\tau/4+\varepsilon}) \label{pos}
\end{eqnarray}

and a similarly 
\begg H_N^*(-):=\sum_{j< 0} \left(G_1^{(j)}+G_2^{(j)}+\ldots
+G_{\xi_2(j,V_N)}^{(j)}\right) 
=\sum_{j< 0}\xi_2(j,V_N)\frac{1-2p_j}{2p_j}+O(N^{3/4-\tau/4+\varepsilon})\label{neg}
\endd

As for the case $\gamma_2=1$ we don't have (\ref{neg}), we need a different argument, as follows.

Recall (\ref{heyde0})  and that in this condition  $1/2< \tau \leq 1.$ 
Select $\delta >0$ such that $\tau=1/2 +2\delta$ should hold.
We will show that for $N$ big enough
\begg H_N^*(-)\leq N^{1-\tau/2+\delta} \quad {\rm a.s.}      \label{delta}
\endd

To this end observe that for N big enough by Lemma C and Lemma E 

$$M _N\leq  N^{1/2} \log N \quad  {\rm a.s.} \quad  {\rm and} \quad \xi(j,V_N)\leq \xi(N) \leq N^{1/2 } \log N \quad  {\rm a.s.}$$ 

where $M_N=max_{ 0 \leq k \leq N} |S_2(k)|.$  Introduce the notation $r_N=N^{1/2} \log N$

implying that

\begg H_N^*(-)\leq \sum_{j=1}^{r_N}
\sum_{i=1}^ {r_N} G^{(-j)}_i .  \quad a.s.\endd

Let  $N_k= k^{k},\quad  \lambda_N=N^{1-\tau/2 +\delta}$  with some $\delta>0.$  Then  $\frac{N_{k+1}}{N_k}\sim e(k+1).$ From Markov inequality  and  (\ref{heyde0})
we have that

\begin{eqnarray}
{\bf P}\left(\sum_{j=1}^{r_{N_{K+1} }}
\sum_{i=1}^ { r_{N_{K+1}} }G^{(-j)}_i >\lambda_{N_K} \right)&\leq&
\frac{{\bf E} \left (
\sum_{j=1}^{   r_{N_{K+1}}}
\sum_{i=1}^ { r_{N_{K+1}}}G^{(-j)}_i \right)}{\lambda_{N_K}}
=\frac{ \left (
\sum_{i=1}^{   r_{N_{K+1}}}
\sum_{j=1}^ { r_{N_{K+1}}}{\bf E}G^{(-j)}_i \right)}
{\lambda_{N_K}} \nonumber \\
&=&\frac{ \sum_{i=1}^{   r_{N_{K+1}}} \sum_{j=1}^ {  r_{N_{K+1} } }  \frac{1-2p_{-j}}{2p_{-j}} }
{\lambda_{N_{K}}}
=\frac{\sum_{i=1}^{   r_{N_{K+1}}}\left(\sum_{j=1}^ {  r_{N_{K+1} } }  \frac{1}{2p_{-j}} -  r_{N_{K+1} } \right)}{\lambda_{N_{K}}}
\nonumber \\
&=& \frac{r_{N_{K+1}}  o(r_{N_{K+1}})^{1-\tau }}{\lambda_{N_{K}}}\leq \frac {(N_{K+1} )^{1-\tau/2+\delta/2 }}
 {(N_K)^{1-\tau/2+\delta}}\nonumber \\
 &\sim& (e(K+1))^{1-\tau/2+\delta}\frac{1}{(K+1)^{\frac{(K+1)\delta}{2}}}.
\end{eqnarray}

So we got  the $(K+1)$-th  term of a convergent series. By Borel-Cantelli and the monotonicity of  $H_N^*(-)$  we conclude that for $N_K\leq N \leq N_{K+1}$

\begin{eqnarray}
H_N^*(-)\leq H_{N_{K+1}}^*(-)\leq \lambda_{N_K}\leq \lambda_{N} =N^{1-\tau/2+\delta} \quad  a.s.
\end{eqnarray}

for N big enough, proving  (\ref{delta}). 
Applying now Lemma 2.1 with  $\rho=1$ and $\beta_j=2p_{-j}, \, j=1,2,\dots$ imply that
$$\sum_{j<0}\xi_2(j,V_N)\frac{1}{2p_j} =V_N(-)+O(N^{5/4-\tau/2+\ep})$$ 
or  equivalently
$$\sum_{j<0}\xi_2(j,V_N)\frac{1-2p_j}{2p_j} =O(N^{5/4-\tau/2+\ep})$$ 
Consequently
\begin{eqnarray}
H_N^*(-)&=&
\sum_{j<0}\xi_2(j,V_N)\frac{1-2p_j}{2p_j} +O(N^{5/4-\tau/2+\ep})+N^{1-\tau/2+\delta} \nonumber \\
&=&
\sum_{j<0}\xi_2(j,V_N)\frac{1-2p_j}{2p_j} +O(N^{5/4-\tau/2+\ep}). \label{secondneg}
\end{eqnarray}
 as we can select $\delta$ to be arbitrary small.  Consequently,  we  have by (\ref{pos}), (\ref{neg}) and (\ref{secondneg}) that for $1\leq\lambda_1<\lambda_2$

\begin{eqnarray}
H_N^*&=&\sum_j \left(G_1^{(j)}+G_2^{(j)}+\ldots
+G_{\xi_2(j,V_N)}^{(j)}\right) \nonumber
\\
&=&\sum_{j}\xi_2(j,V_N)\frac{1-2p_j}{2p_j}+O(N^{3/4-\tau/4+\varepsilon})+O(N^{5/4-\tau/2+\ep}).
\nonumber
 \\&=&-V_N+\sum_j \xi_2(j,V_N)\frac{1}{2p_j}
+O(N^{3/4-\tau/4+\varepsilon})++O(N^{5/4-\tau/2+\ep}).\nonumber
 \\&=&-V_N+ \widehat{A}_2(V_N)+O(N^{3/4-\tau/4+\varepsilon})+O(N^{5/4-\tau/2+\varepsilon}).\nonumber
 \\&=&-V_N+A_2(V_N)
+O(N^{5/4-\tau/2+ \varepsilon}) +O(N^{3/4+\ep}) \nonumber
\\ &=& -V_N+A_2(V_N)
+O(N^{5/4-\tau/2+ \varepsilon}) \quad a.s. \label{fontos}
\end{eqnarray}
where we used the definition of   $\widehat{A}_2(.)$  and Lemma J.

Clearly, using   (\ref{haen}) and (\ref{fontos})

 $$N=H_N+V_N=H_N^*+V_N+O(N^{1/2+\delta})=A_2(V_N)
+O(N^{5/4-\tau/2+\varepsilon}) \quad a.s.$$

and
$$V_N=A_2^{-1}(N)+O(N^{5/4-\tau/2+\ep}) \qquad  a.s.$$

\bigskip
\noindent
{\bf Remark 3.1} In the previous line we used the fact that 
$A_2^{-1}(u+v)-A_2^{-1}(u)\leq v.$  
To see this, first recall from Lemma 3.1 that $A_2(t),\,$  $A_2^{-1} (t)$  and  
$\alpha(t) =A_2(t)-t$ are all nondecreasing. Then
$$v=A_2(A_2^{-1}(u+v))-A_2(A_2^{-1}(u))=\alpha(A_2^{-1}(u+v))+A_2^{-1}(u+v)
-\alpha(A_2^{-1}(u))-A_2^{-1}(u)$$
$$\geq A_2^{-1}(u+v)-A_2^{-1}(u).$$

So we can conclude, using Lemmas H and I that 
$$C_2(N)=S_2(V_N)=W_2(V_N)+O(N^{1/4+\ep})=W_2((A_2^{-1}(N)
+O(N^{5/4-\tau/2+\ep}))+O(N^{1/4+\ep})$$
$$=W_2((A_2^{-1}(N))+O(N^{5/8-\tau/4+\ep}) \qquad  a.s.$$
$$C_1(N)=S_1(H_N)=S_1(N-V_N)=W_1(N-V_N)+O(N^{1/4+\ep})
=W_1(N-A_2^{-1}(N))+O(N^{5/8-\tau/4+\ep})) \qquad  a.s.$$
proving our theorem. $\Box$

\bigskip
\noindent
{\bf Remark 3.2}  We could extend our result for  the case $\gamma_1=\gamma_2=1$, by using the argument  
of the case $\gamma_2=1$ for both of the positive and the negative side. Then we would get the following result:

{\it Under the conditions {\rm (\ref{heyde0})}, and   
$\gamma_1 = \gamma_2=1$
on an appropriate probability space for the random walk $\{{\bf
C}(N)=(C_1(N),C_2(N));\, \, N=0,1,2,\ldots\}$ one can construct two independent
standard Wiener processes $\{W_1(t);\, t\geq 0\}$, $\{W_2(t);\,
t\geq 0\}$ so that, as $N\to\infty$, we have with any}
$\varepsilon>0$
\begg
|C_1(N)|+|C_2(N)-W_2(N)|
=O(N^{5/8-\tau/4+\varepsilon})\quad a.s. \label{hat}
\endd

\noindent
{\bf Remark 3.3} In our paper (\cite{CF21}  Lemma 4.1) we calculated 
 the density function of  $A^{-1}(t)$   and $t-A^{-1}(t).$ Our result was the following:

\bigskip

{\it  Suppose that }$\gamma_1>\gamma_2\geq 1.$ 
$$P(A^{-1}(t)\in dv) =\frac{t}{\pi v}\frac{1}{\sqrt{(v\gamma_1-t)(t-\gamma_2v)}}\,dv\quad {\rm for}\quad \frac{t}{\gamma_1}<v<\frac{t}{\gamma_2},$$
\begin{eqnarray}P(t-A^{-1}(t)\in dv) =
\frac{t}{\pi(t-v)}\frac{1}{\sqrt{((\gamma_1-1)t-\gamma_1 v)(t(1-\gamma_2) +\gamma_2v)}}\,dv  \nonumber\\
\quad {\rm for}\quad t \left(1-\frac{1}{\gamma_2}\right)<v<t\left(1-\frac{1}{\gamma_1}\right).\nonumber
\end{eqnarray}

 As in  our Theorem 1.1 our random walk $\{{\bf
C}(N)=(C_1(N),C_2(N)$ is approximated  with the same pair of oscillating Wiener processes as in Theorem D, 
 we get the same consequences as in case of the $p_i$-s were restricted to be 1/2 or 1/4.
We proved the following laws of the iterated logarithm  (\cite{CF21}). 

\bigskip
\noindent
{\bf Corollary A} (\cite{CF21}  Corollary 4.1).
{ \it Under condition {\rm (\ref{heyde0} )} with   $\gamma_1>\gamma_2\geq 1$ the following laws of the iterated logarithm hold.}
\begin{noitemize}
\item
$$
\limsup_{t\to\infty}\frac{W_1(t-A^{-1}(t))}{\sqrt{t\log\log t}}=
\limsup_{N\to\infty}\frac{C_1(N)}{\sqrt{N\log\log N}}=\sqrt{2\left(1-\frac{1}{\gamma_1}\right)}
\quad a.s.,$$
\item
$$\liminf_{t\to\infty}\frac{W_1(t-A^{-1}(t))}{\sqrt{t\log\log t}}=
\liminf_{N\to\infty}\frac{C_1(N)}{\sqrt{N\log\log N}}=-\sqrt{2\left(1-\frac{1}{\gamma_1}\right)}
\quad a.s.,$$
\item
$$
\limsup_{t\to\infty}\frac{W_2(A^{-1}(t))}{\sqrt{t\log\log t}}=
\limsup_{N\to\infty}\frac{C_2(N)}{\sqrt{N\log\log N}}=\sqrt{\frac{2}{\gamma_1}}
\quad a.s.,
$$
\item
$$
\liminf_{t\to \infty}\frac{W_2(A^{-1}(t))}{\sqrt{t\log\log t}}=
\liminf_{N\to\infty}\frac{C_2(N)}{\sqrt{N\log\log N}}=-\sqrt{\frac{2}{\gamma_2}}
\quad a.s.
$$
\end{noitemize}


\bigskip


 
 
\begin{thebibliography}{9}
\bibitem {BE} \textsc{Bertacchi, D.} (2006). Asymptotic behavior of the
simple random walk on the 2-dimensional comb. \textit{Electron. J.
Probab.} \textbf{11} 1184--1203.

\bibitem{BZ}
\textsc{Bertacchi, D. and Zucca, F.} (2003). Uniform asymptotic
estimates of transition probabilities on combs. \textit{J. Aust. Math.
Soc.} \textbf{75} 325--353.

\bibitem{CX}
\textsc{Chen, X.} (1999). How often does a Harris recurrent Markov chain
recur? \textit{Ann. Probab.} \textbf{27} 1324--1346.

\bibitem{C00}
\textsc{Chen, X.} (2000). On the limit laws of second order for additive
functionals of Harris recurrent Markov chain. \textit{Probab. Th. Rel.
Fields} \textbf{116} 89--123.

\bibitem{C01}
\textsc{Chen, X.} (2001). Moderate deviations for Markovian occupation
time. \textit{Stochastic Process. Appl.} \textbf{94} 51--70.

\bibitem{CH}
\textsc{Chung, K.L.} (1948). On the maximum partial sums of
sequences of independent random variables. \textit{Trans. Amer.
Math. Soc.} \textbf{64} 205--233.

\bibitem{CO}
\textsc{Comtet, L.} (1974). \textit{Advanced Combinatorics. The Art of
Finite and Infinite Expansions}, D. Reidel Publishing Co., Dordrecht,
enlarged edition.

\bibitem{CCFR08}
\textsc{Cs\'aki, E., Cs\"org\H o, M., F\"oldes, A. and R\'ev\'esz,
P.} (2009). Strong limit theorems for a simple random walk on the
2-dimensional comb. \textit{Electr. J. Probab.} \textbf{14} 2371--2390.

\bibitem{CCFR12}
\textsc{Cs\'aki, E., Cs\"org\H o, M., F\"oldes, A. and R\'ev\'esz, P.} (2012).
Random walk on half-plane half-comb structure. \textit{Ann. Math. Inform.}
\textbf{39} 29--44. MR2959879

\bibitem{CCFR13} 
\textsc{Cs\'aki, E., Cs\"org\H o, M., F\"oldes, A. and R\'ev\'esz, P.} (2013). 
Strong limit theorems for anisotropic random walks on $Z^2$.  
\textit{Periodica Math. Hungar.} \textbf{67} 71-94. 

\bibitem{CF21} \textsc{Cs\'aki, E. and F\"oldes, A.} (2021)
Random walks on Comb-Type Subsets of $\mathbb{Z}^2$ 
Journal of Theoretical Probability) 1-25


\bibitem{CR83}\textsc{Cs\'aki, E. and R\'ev\'esz, P.} (1983). Strong
invariancefor local time. \textit{Z. Wahrsch. verw. Gebiete} \textbf{50}
5--25.


\bibitem{CRR}
\textsc{Cs\'aki, E., R\'ev\'esz, P. and Rosen, J.} (1998).
Functional laws of the iterated logarithm for local times of recurrent
random walks on $\mathbb{Z}^2$. \textit{Ann. Inst. H. Poincar\'e,
Probab. Statist.} \textbf{34} 545--563.


\bibitem{CSR79}
\textsc{Cs\"org\H o, M. and R\'ev\'esz, P.} (1979) How big are the increments of a 
Wiener process? \textit{Ann. Probab.} \textbf{7} 731-737.
\bibitem{DK}
\textsc{Darling, D.A. and Kac, M.} (1957). On occupation times for
Markoff processes. \textit{Trans. Amer. Math. Soc.} \textbf{84}
444--458.

\bibitem{DH}
\textsc{den Hollander, F.} (1994). On three conjectures by K. Shuler.
\textit{J. Statist. Physics} \textbf{75} 891--918.

\bibitem{DE}
\textsc{Dvoretzky, A. and Erd\H os, P.} (1951). Some problems on random
walk in space. \textit{Proc. Second Berkeley Symposium}, pp. 353--367.

\bibitem{ET}
\textsc{Erd\H os, P. and Taylor, S.J.} (1960). Some problems concerning
the structure of random walk paths. \textit{Acta Math. Acad. Sci.
Hungar.} \textbf{11} 137--162.

\bibitem{H}
\textsc{Heyde, C.C.} (1982). On the asymptotic behaviour of random walks
on an anisotropic lattice. \textit{J. Statist. Physics} \textbf{27}
721--730.

\bibitem{H93}
\textsc{Heyde, C.C.} (1993). Asymptotics for two-dimensional anisotropic
random walks. In: \textit{Stochastic Processes.} Springer, New York,
pp. 125--130.

\bibitem{HWW}
\textsc{Heyde, C.C., Westcott, M. and Williams, E.R.} (1982).
The asymptotic behavior of a random walk on a dual-medium lattice.
\textit{J. Statist. Physics} \textbf{28} 375-380.

\bibitem{K}\textsc{Kesten, H.} (1965) An iterated logarithm law for
the local time. \textit{Duke Math.J.} \textbf{32} 447-456.


\bibitem{KMT}
\textsc{Koml\'os, J., Major, P. and Tusn\'ady, G.} (1975).
An approximation of partial sums of independent rv's and the sample df.
I. \textit{Z. Wahrsch. verw. Gebiete} \textbf{32} 111--131.

\bibitem{RE81}
\textsc{R\'ev\'esz, P.} (1981). Local time and invariance. 
\textit{Lecture Notes in Math.} \textbf{861} 128--145. Springer, New York.


\bibitem{RE}
\textsc{R\'ev\'esz, P.} (2005).
\textit{Random Walk in Random and Non-Random Environments}, 2nd ed.
World Scientific, Singapore.



\bibitem{SLS}
\textsc{Seshadri, V., Lindenberg, K. and Shuler, K.E.} (1979).
Random Walks on Periodic and Random Lattices. II. Random Walk Properties
via Generating Function Techniques. \textit{J. Statist. Physics}
\textbf{21} 517--548.

\bibitem{SH}
\textsc{Shuler, K.E.} (1979). Random walks on sparsely periodic and random
lattices I. \textit{Physica A} \textbf{95} 12--34.

\bibitem{SSL}
\textsc{Silver, H., Shuler, K.E. and Lindenberg, K.} (1977).
Two-dimensional anisotropic random walks. In:
\textit{ Statistical mechanics and statistical methods in theory and
application (Proc. Sympos., Univ. Rochester, Rochester, N.Y., 1976)},
Plenum, New York, pp. 463--505.

\bibitem{WE}
\textsc{Westcott, M.} (1982). Random walks on a lattice. \textit{J.
Statist. Physics} \textbf{27} 75--82



\bibitem{WH}
\textsc{Weiss, G.H. and Havlin, S.} (1986). Some properties of a random
walk on a comb structure. \textit{Physica A} \textbf{134} 474--482.
\end{thebibliography}
\end{document}